\title{Tiling transitive tournaments and their blow-ups}

\author{Raphael Yuster
\thanks{e-mail: raphy@research.haifa.ac.il \qquad
World Wide Web: http:$\backslash\backslash $research.haifa.ac.il$\backslash$\symbol{126}raphy}
\\ Department of Mathematics\\ University of
Haifa at Oranim\\ Tivon 36006, Israel}

\date{} 

\documentstyle [11pt]{article}
\newtheorem{theo}{Theorem}[section]
\newtheorem{prop}[theo]{Proposition}
\newtheorem{lemma}[theo]{Lemma}
\newtheorem{coro}[theo]{Corollary}

\newcommand\npf{\mbox{ }\hfill\sqr\vskip6pt}
\def\sqr{$\vcenter{\hrule height.2mm
\hbox{\vrule width.2mm height2mm\kern2mm
\vrule width.2mm}\hrule height.2mm}$}
\newcommand{\ignore}[1]{}

\begin{document}
\maketitle
\setcounter{page}{1}
\begin{abstract}
Let $TT_k$ denote the transitive tournament on $k$ vertices.
Let $TT(h,k)$ denote the graph obtained from $TT_k$ by replacing each vertex with an independent set
of size $h \geq 1$. The following result is proved:
Let $c_2=1/2$, $c_3=5/6$ and $c_k=1-2^{-k-\log k}$ for $k \geq 4$.
For every $\epsilon > 0$ there exists $N=N(\epsilon,h,k)$ such that for every undirected graph $G$ with
$n > N$ vertices and with $\delta(G) \geq c_kn$, every orientation of $G$
contains vertex disjoint copies of $TT(h,k)$ that cover all but at most $\epsilon n$ vertices.
In the cases $k=2$ and $k=3$ the result is asymptotically tight.
For $k \geq 4$, $c_k$ cannot be improved to less than $1-2^{-0.5k(1+o(1))}$.
\end{abstract}

\section{Introduction}
All graphs considered here are finite and simple. For standard
terminology on undirected and directed graphs the reader is
referred to \cite{Bo}. Finding many isomorphic copies of a given graph $H$ within a larger graph
$G$ is a central topic in extremal graph theory that has been studied extensively
in recent years. Formally, a graph $G$ has an {\em $H$-factor} if it contains
a spanning subgraph whose components are isomorphic to $H$.
For $ 0 \leq \alpha \leq 1$ we say that a graph $G$ of order $n$ has an
{\em $(H,\alpha)$-factor} if there are vertex disjoint copies of $H$ in $G$ that cover $\alpha n$
vertices of $G$. Thus, an $(H,1)$-factor is an $H$-factor. Most of the results
on $H$-factors and almost $H$-factors (namely, results guaranteeing $(H,1-\epsilon)$-factors)
are stated in terms of the chromatic number of $H$, or closely related variants of the chromatic number.
Perhaps the most important result is that of Hajnal and Szemer\'edi \cite{HaSz} 
stating that an $n$-vertex graph with minimum degree at least $n(1-1/k)$ has a $K_k$-factor,
assuming $k$ divides $n$.
This result was extended by Alon and Yuster in \cite{AlYu2} to arbitrary graphs $H$ with $\chi(H)=k$
at the price of increasing the minimum degree requirement by $\epsilon n$, and having $n$ sufficiently
large. Later, Koml\'os, S\'ark\"ozi and Szemer\'edi \cite{KoSaSz} showed that $\epsilon n$ can be replaced with
a constant depending only on $H$. Koml\'os proved an almost $H$-factor result
which is stated in terms of the {\em critical chromatic number} of $H$. This parameter,
which is greater than $\chi(H)-1$ and is at most $\chi(H)$,
takes into account the fact that the vertex class sizes of an optimal coloring may vary significantly.
His result was extended recently by Shokoufandeh and Zhao \cite{ShZa}.

Let $H$ be a digraph. In this paper we study the $H$-factor and almost $H$-factor problems in orientations
of an undirected graph $G$. We say that $H$ is {\em immuned} against orientations of
$G$ if every orientation of $G$ contains $H$. Clearly the definition is interesting only if $H$ is acyclic.
This problem has been extensively studied when $G$ is a complete graph.
Let $TT_k$ denote the unique transitive tournament with $k$ vertices.
Let $f(k)$ denote the minimum integer $n$ that guarantees that every orientation of $K_n$
has a $TT_k$. A trivial induction argument gives $f(k) \leq 2^{k-1}$.
On the other hand, Erd\H{o}s and Moser \cite{ErMo} proved, using the probabilistic method, that
$f(k) \geq 2^{0.5k(1+o(1))}$. It is easy to show $f(2)=2$, $f(3)=4$ an it is well known
that and $f(4)=8$ and $f(5)=14$ \cite{RePa}. It is currently known that $f(7) \leq 54$
\cite{Sa} and, therefore, the induction argument gives $f(k) \leq 54 \cdot 2^{k-7}$.
Let $g(k)$ denote the minimum integer $n$ that guarantees that every orientation of $K_n$
has a $TT_k$-factor, assuming $k$ divides $n$.
The (nontrivial) existence of $g(k)$ is attributed to Erd\H{o}s in \cite{Re} and
a (huge) upper bound yielding $g(k) < 2^{k2^k}$ follows from Lonc and Truszcy\'nski \cite{LoTr}.
This upper bound was significantly improved to $g(k) < 4^k$ by Caro \cite{Ca}, but is
probably still far from being optimal. Chen, Lu and West \cite{ChLuWe} proved that every orientation of
$K_{4m^2-6m}$ has an $H$-factor where $H$ is the star with $m$ vertices where all
$m-1$ edges either all emanate from or enter the root.

As in the undirected case, if $G$ is not necessarily complete, the existence of an almost $H$-factor
or an $H$-factor in every orientation of $G$ may be guaranteed only if $G$ has a sufficiently high minimum degree.
However, the discussion in the previous paragraph suggests that the required minimum degree is much
larger than in the undirected case. Moreover, unlike the undirected case, it is impossible to state a
minimum degree condition in terms of the chromatic number alone. To see this, consider even the simplest
case where $H=K_{1,m}$ is the star with $m+1$ vertices and all edges emanate from the root.
Let $G$ be the complete $m+1$-partite graph with all vertex classes having the same size $h$,
except the first one that has $h(1+\alpha)$ vertices and the second one that has $h(1-\alpha)$ vertices.
It is easy to see that $n=(m+1)h$ is the number of vertices of $G$, $\delta(G) \geq (m-\alpha)h$.
However, the orientation of $G$ in which all edges are directed from lower indexed classes to higher ones
has at least $\alpha h$ uncovered vertices (all from the first class) in any set of vertex disjoint copies of $H$.
Thus, although $\chi(H)=2$, the minimum degree of $G$ can be arbitrarily large for $m$ sufficiently large.
Fortunately, for some important classes of digraphs, the chromatic number implies a minimum degree bound for the
existence of an almost $H$-factor. We say that an acyclic digraph $H$ has a {\em balanced $k$-coloring}
if it can be properly colored with $k$ colors such that: (i) every color class has the same number of vertices
and (ii) all the edges between any two color classes are in the same direction.
Denote by ${\vec \chi(H)}$ the minimum number of colors in a balanced coloring of $H$ (in case one exists).
For example, ${\vec \chi(TT_k)}=k$. Similarly, ${\vec \chi(TT(h,k))}=k$ where $TT(h,k)$ is the
complete $k$-partite acyclic digraph with $h$ vertices in each part, and with all edges between and two parts going in the
same direction. Another example: the unique graph $H$ obtained by orienting a path on three edges such that
there is no directed path of length 2 has ${\vec \chi(H)}=2$. Our main result is the following:
\begin{theo}
\label{t1}
Let $c_2=1/2$, $c_3=5/6$ and $c_k=1-2^{-k-\log k}$ for $k \geq 4$.
For every digraph $H$ having ${\vec \chi(H)}=k$
and for every $\epsilon > 0$,
there exists $N=N(\epsilon,H)$ such that for every undirected graph $G$ with
$n > N$ vertices and with $\delta(G) \geq c_kn$, every orientation of $G$
contains vertex disjoint copies of $H$ that cover all but at most $\epsilon n$ vertices.
\end{theo}
Notice that it suffices to prove Theorem \ref{t1} for the graphs $TT(h,k)$.
The tightness of the result for $k=2$ is trivial. For any $\gamma > 0$ the complete bipartite graph $G$
with $(1/2-\gamma)n$ vertices in one vertex class and $(1/2+\gamma)n$ vertices in the other vertex class
does not contain more than $(1/2-\gamma)n$ vertex-disjoint copies of, say, $TT_2$, in any orientation of
$G$. Hence we always remain with at least $2\gamma n$ vertices uncovered.
The proof of tightness of the $k=3$ case is slightly more complicated.
\begin{prop}
\label{p2}
For every $1/30 > \gamma > 0$ there exists an $n$-vertex graph $G$ with $\delta(G) \geq (5/6-\gamma)n$ and
an orientation of $G$ having at most $n/3-\gamma n$ vertex-disjoint copies of $TT_3$.
\end{prop}
For larger $k$, the constant $c_k=1-2^{-k-\log k}$ cannot be dramatically improved, in the sense that
one cannot replace it with a constant smaller than $1-2^{-0.5k(1+o(1))}$.
This follows easily from the above mentioned lower bound for $f(k)$.
We elaborate more on this fact in Section 4.

We now turn to the problem of finding an exact $H$-factor. Extending Caro's proof mentioned above,
stating that orientations of $K_n$ have a $TT_k$-factor whenever $n > 4^k$ is a multiple of $k$,
we can prove the following ``dense graph'' version.
\begin{theo}
\label{t2}
Let $H$ be an acyclic digraph with $h$ vertices. If $G$ has $n$ vertices,
$\delta(G) \geq n(1-1/4^h)+4^h$, and $h | n$, then every orientation of $G$ has an $H$-factor.
\end{theo}

The rest of this paper is organized as follows. In Section 2 we present the necessary tools for the proof
of Theorem \ref{t1}. Section 3 contains the proof of Theorem \ref{t1}. Section 4 considers the lower
bounds for $c_k$ and the proof of Proposition \ref{p2}. Section 5 considers exact $H$-factors and contains
the proof of Theorem \ref{t2}.
The final section contains some concluding remarks and open problems.

\section{Lemmas and tools}
Let $K(t,r)$ denote the complete $r$-partite graph with $t$ vertices in each partite class.
In the proof of Theorem \ref{t1} it will be useful to show that for $r=r(k)$ that is relatively small,
$K(t,r)$ contains ``many'' vertex-disjoint copies of $TT_k$ in any orientation of the edges of $K(t,r)$.
By ``many'' we mean that the number of uncovered vertices is independent of $t$.
This is trivial for $k=2$ since any orientation of $K(t,2)$ trivially has a $TT_2$-factor.
It is also easy for $k=3$. Every orientation of $K_6=K(1,6)$ is easily verified to contain two vertex disjoint copies of
$TT_3$. Thus, an oriented $K(t,6)$ has a $TT_3$-factor.
In fact, an immediate consequence of the proof of Proposition \ref{p2} is that for every constant $C$,
for $t$ sufficiently large, there are orientations of $K(t,5)$ such that in  every maximal set of vertex-disjoint
$TT_3$ there remain at least $C$ uncovered vertices. Hence the choice $r=6$ is best possible for $k=3$.
If we wish to guarantee no loss at all, that is, a $TT_k$-factor, then $r$ would grow too large.
The best known value for $r$ in this case would be as large as Caro's  upper bound for $g(k)$ mentioned
in the introduction, and which is close to $4^k$. If we settle for an {\em almost} factor we can do much better.
Let $f^*(k)$ denote the minimum integer $m$ that guarantees that in any orientation of $K_m$, and for
every vertex of $K_m$, there is a $TT_k$ containing the vertex. Recalling the definition of $f(k)$ mentioned
in the introduction, we clearly have $f^*(k) \geq f(k)$. An easy inductive argument yields $f^*(k) \leq 2^{k-1}$.
\begin{lemma}
\label{l21}
Let $r_2=2$, $r_3=6$ and $r_k=k(f^*(k)-2)+2$ for $k \geq 4$.
In any orientation of $K(t,r_k)$ there are vertex-disjoint copies of $TT_k$ that cover all
but at most $f^*(k)-1$ vertices.
\end{lemma}
{\bf Proof}\,
As shown above, we only need to prove the lemma for $k \geq 4$.
We prove something slightly stronger. In any orientation of $K(t,r_k)$ there are vertex-disjoint copies of $TT_k$ that cover all
but at most $f^*(k)-1$ vertices, and the uncovered vertices induce a complete graph.
We use induction on $t$.  The case $t=1$ is trivial from the definition of $f^*(k)$.
Assuming the lemma holds for $K(t-1,r_k)$, we prove it for $K(t,r_k)$. Fix an orientation of
$K(t,r_k)$. Delete one vertex from each partite class, and find in the resulting $K(t-1,r_k)$
a set of vertex-disjoint copies of $TT_k$ satisfying the assertion. There are at most $f^*(k)-1$
uncovered vertices in the $K(t-1,r_k)$, each belonging to a distinct partite class, and we also have
the $r_k$ uncovered originally deleted vertices. For each uncovered vertex of the $K(t-1,r_k)$ we pick
a copy of $TT_k$ containing it, and $k-1$ of the originally deleted vertices. This can be done even
for the last uncovered vertex of the $K(t-1,r_k)$ since up till now we only used at most
$(k-1)(f^*(k)-2)$ originally deleted vertices and we therefore still have at least
$$
r_k-(k-1)(f^*(k)-2) =k(f^*(k)-2)+2 - (k-1)(f^*(k)-2) = f^*(k)
$$
deleted vertices in our disposal. Now, we may remain with several uncovered originally deleted
vertices. As long as there are at least $f^*(k)$ of them, we can greedily select another $TT_k$.
After the end of the process we remain with at most $f^*(k)-1$ uncovered vertices that belong
each to a distinct partite class. \npf

\noindent
Lemma \ref{l21} is not optimal. For $k=4$ we have $f^*(4)=8$ and the lemma
gives a bound of $r_4=26$. A tailor-made proof for the case $k=4$ using the same arguments
works already for $r_4=20$. (Since both 20 and $f^*(4)$ are multiples of $4$ we can assume in the proof that we
always remain with at most $f^*(k)-k$ uncovered vertices in this case.)

An important tool used in the proof of Theorem \ref{t1} is the following directed version of Szemer\'edi's regularity
lemma. Although never published, this lemma is a relatively easy consequence
of the standard regularity lemma proved in \cite{Sz}.
For more details on the regularity lemma we refer the
reader to the excellent survey of Koml\'os and Simonovits \cite{KoSi}, which
discusses various applications of this powerful result,
and to \cite{AlYu1} which addresses another problem solved with the aid of the directed version of the lemma.
We now give the definitions necessary in order to state the directed regularity lemma.

Let $G=(V,E)$ be a directed graph, and let $A$ and $B$ be two disjoint subsets
of $V(G)$. If $A$ and $B$ are non-empty and $e(A,B)$ is the number of edges from $A$ to $B$,
define the {\em density of edges from $A$ to $B$} as
$$
d(A,B) = \frac{e(A,B)}{|A||B|}.
$$
For $\gamma >0$ the pair $(A,B)$ is called {\em $\gamma$-regular}
if for every $X \subset A$ and $Y \subset B$ satisfying
$|X|>\gamma |A|$ and $|Y|>\gamma |B|$ we have
$$
|d(X,Y)-d(A,B)| < \gamma \qquad \qquad |d(Y,X)-d(B,A)| < \gamma.
$$

An {\em equitable partition} of a set $V$ is a partition of $V$ into
pairwise disjoint classes $V_1,\ldots,V_m$ whose sizes are as equal as possible.
An equitable partition
of the set of vertices $V$ of a directed graph $G$ into the classes $V_1,\ldots,V_m$ is
called {\em $\gamma$-regular} if $|V_i| \leq \gamma |V|$ for
every $i$ and all but at most $\gamma {m \choose 2}$ of the pairs
$(V_i,V_j)$ are $\gamma$-regular.

The directed regularity lemma states the following:
\begin{lemma}
\label{l22}
For every $\gamma>0$, there is an integer $M(\gamma)>0$ such that for
every directed graph $G$ of order $n > M$ there is a $\gamma$-regular partition of
the vertex set of $G$ into $m$ classes, for some $1/\gamma \leq m \leq M$. \npf
\end{lemma}

A useful notion associated with a $\gamma$-regular partition is that of the {\em {\bf or} cluster graph}.
Suppose that $G$ is a directed graph with a $\gamma$-regular partition
$V= V_1 \cup \cdots \cup V_m$, and $\eta>0$
is some fixed constant (to be thought of as small, but much larger than $\gamma$).
The {\em undirected} {\bf or} cluster graph $C(\eta)$ is defined on the vertex set $\{1,\ldots,m\}$ by
declaring $ij$ to be an edge if
$(V_i,V_j)$ is a $\gamma$-regular pair with $d(V_i,V_j) \geq \eta$ {\em or}
$d(V_j,V_i) \geq \eta$ (in some applications such as that appearing in \cite{AlYu1} one needs to use
the analogous {\em {\bf and} cluster graph}).

Our next tool is the following result of Alon and Yuster \cite{AlYu2}, extending the theorem of Hajnal and Szemer\'edi
\cite{HaSz} to complete partite graphs.
\begin{lemma}
\label{l23}
Let $t$ and $r$ be positive integers and let $\beta > 0$. There exists $M^*=M^*(t,r,\beta)$
such that every undirected graph $C^*$ with $m^* > M^*$ vertices where $tr$ divides $m^*$ and
$\delta(C^*) \geq m^*(1-1/r+\beta)$ has a $K(t,r)$-factor. \npf
\end{lemma}
We note here that the proof of Lemma \ref{l23} also uses the (undirected) regularity lemma.
In the proof of Theorem \ref{t1} we apply Lemma \ref{l23} to a subgraph of the cluster
graph resulting from the application of the directed regularity lemma. Hence,
the proof of Theorem \ref{t1} requires, essentially, a double application of the regularity lemma.

We shall require the following corollary of Lemma \ref{l23}.
\begin{coro}
\label{c24}
Let $t$ and $r$ be positive integers and let $1/(5r) > \beta > 0$. There exists $T=T(t,r,\beta)$
such that the following holds. If $C$ is an undirected graph with $m > T$
vertices in which the degrees of all vertices but at most $\beta m$ are at least
$m(1-1/r-\beta)$ then $C$ contains a set of at least $\frac{m}{tr}-6\beta r m$ vertex disjoint
copies of $K(t,r)$.
\end{coro}
{\bf Proof}\,
Put $T = \max\{t/\beta, M^*(t,r,\beta)\}$ where $M^*$ is the constant from Lemma \ref{l23}.
Let $C$ be a graph with $m > T$ vertices satisfying the conditions of the corollary.
Let $V'$ be the set of all vertices of $C$ whose degrees in $C$
are less than $m(1-1/r-\beta)$. Let $C'$ be the graph 
obtained from $C$ by joining each vertex of $V'$ to any other vertex
of $C$. (Thus in $C'$ the degree of each vertex in $V'$ is $m-1$).
Let $C^*$ be the graph obtained from $C'$ by adding to it
a complete graph on a set $V^*$ of at least $4\beta rm$ and at most $5\beta rm$
new vertices and by joining each of them to every vertex of $C'$.
The exact size of $V^*$ is chosen so that the total number
of vertices of $C^*$ will be divisible by $tr$. In $C^*$ the degree of 
every vertex in $V' \cup V^*$ is $m^*-1$, where $m^*=m+|V^*|$
is the number of vertices of $C^*$. The degree of each other vertex
is at least $(1-1/r-\beta)m+|V^*|$ and notice that
$$
\left(1-\frac{1}{r}-\beta\right)m+|V^*| =
\left(1-\frac{1}{r}-\beta\right)m^*+(m^*-m)\left(\frac{1}{r}+\beta\right) \geq \left(1-\frac{1}{r}-\beta\right)m^*+4\beta m
$$
$$
\geq \left(1-\frac{1}{r}-\beta\right)m^*+2\beta m^* = \left(1-\frac{1}{r}+\beta\right)m^*.
$$
Therefore, by Lemma \ref{l23}, $C^*$ has a set of $m^*/(tr)$ vertex disjoint
copies of $K(t,r)$. At most $|V'|+|V^*| \leq 6\beta r m$ of these copies
contain vertices of $V' \cup V^*$ and all the others are in fact 
subgraphs of $C$. Therefore, $C$ contains a set of at least
$m^*/(tr)-|V'|-|V^*| \geq m/(tr) - 6\beta r m$ vertex disjoint
copies of $K(t,r)$. \npf

\section{Proof of Theorem \ref{t1}}
Recall that $c_2=1/2$, $c_3=5/6$ and $c_k=1-2^{-k-\log k}$ for $k \geq 4$
and recall the definition of $r_k$ from lemma \ref{l21}, where $r_2=2$, $r_3=6$ and
$r_k=k(f^*(k)-2)+2 \leq k(2^{k-1}-2)+2 < 2^{k + \log k}$ for $k \geq 4$.
Thus, $1-1/r_k \leq c_k$ for all $k \geq 2$. Also recall that it suffices to prove
Theorem \ref{t1} for the graphs $TT(h,k)$. Hence, it suffices to prove
the following slightly stronger version of Theorem \ref{t1}.
\begin{theo}
\label{t31}
Let $h \geq 1$ and $k \geq 2$ be positive integers and let $\epsilon > 0$.
There exists $N=N(\epsilon,h.k)$ such that for every undirected graph $G$ with
$n > N$ vertices and with $\delta(G) \geq n(1-1/r_k)$, every orientation of $G$
contains vertex disjoint copies of $TT(h,k)$ that cover all but at most $\epsilon n$ vertices.
\end{theo}
{\bf Proof}\,
We first select constants $\eta, t, \mu, \gamma, M, N$ as follows.
Let $\eta=\epsilon/(200r_k)$. Let $t$ be the smallest integer satisfying $84\eta r_k + 2^k/(tr_k) < \epsilon /2$.
Let $\mu=(\eta/4)^{2hk}$. Let $T(t,r_k,7\eta)$ be the constant from Corollary \ref{c24} and let
$$
\gamma=\min\left\{\mu^2 ~,~\frac{1}{T(t,r_k,7\eta)}\right\}.
$$
Let $M=M(\gamma)$ be defined as in Lemma \ref{l22}. Let $N=\max\{3M^2 ~,~ M/\mu^2\}$.

Let $G$ be an undirected graph with $n > N$ vertices. Fix an orientation ${\vec G}$ of $G$.
We begin by applying Lemma \ref{l22} to ${\vec G}$.
Notice that $n>N > M(\gamma)$ so Lemma \ref{l22} yields
a $\gamma$-regular partition of
the vertex set of $G$ into $m$ classes, where $1/\gamma \leq m \leq M$.
Denote the vertex classes by $V_1,\ldots,V_m$. Fix the associated {\bf or} cluster
graph $C=C(\eta)$ on the vertices $\{1,\ldots,m\}$.

We now show that $C$ has a very large subgraph with high minimum degree.
A vertex $i$ of $C$ is called {\em good} if there are at most $\gamma m$ 
other vertices $j$ of $C$ such that the pair $(V_i,V_j)$ is not
$\gamma$-regular. Obviously, all vertices of $C$ but at most $\gamma m$ are good.

\noindent
{\bf Claim:}\,
The degree of any good vertex of $C$ is at least $(1-1/r_k-7\eta)m$.

\noindent
{\bf Proof}\,
Let $b=\lfloor \frac{n}{m} \rfloor$.  Note that the number of 
vertices in each of the sets $V_i$, $1 \leq i \leq m$ is either $b$ or $b+1$.
For each fixed $i$, $1 \leq i \leq m$, the sum of the degrees
in $G$ of the vertices in $V_i$ is at least $(1-1/r_k)n b$,
by the hypotheses. On the other hand, if the degree of $i$ in
$C$ is $d_i$, and $i$ is a good vertex, then the sum of the degrees in 
$G$ of the vertices in $V_i$ can be bounded by the sum of 
four summands, as described below.
\begin{itemize}
\item
The contribution of edges joining two vertices of $V_i$
does not exceed ${{b+1} \choose 2} < b^2$.
\item
The contribution of edges between $V_i$ and classes $V_j$ for which
the pair $(V_i,V_j)$ is not $\gamma$-regular is at most
$(b+1)^2$ times the number of such indices $j$ and is thus at most
$\gamma m (b+1)^2$. (Here we used the fact that $i$ is a good vertex 
of $C$.)
\item
The contribution of edges between $V_i$ and classes $V_j$ for 
which $d(V_i,V_j) < \eta$ and $d(V_j,V_i) < \eta$ does not exceed
$2 m \eta (b+1)^2$.
\item
The contribution of edges between $V_i$ and classes $V_j$ for which
$(V_i,V_j)$ is $\gamma$-regular and either $d(V_i,V_j) \geq \eta$
or $d(V_j,V_i) \geq \eta$
is at most $d_i(b+1)^2$ (since each such $j$ is a neighbor of $i$ in $C$).
\end{itemize}
Therefore
$$
\left(1-\frac{1}{r_k}\right)nb \leq b^2+ \gamma m (b+1)^2 +
2 m \eta (b+1)^2 + d_i(b+1)^2 < b^2+2\gamma m b^2+3m\eta b^2+d_ib(b+3).
$$
Since $b \leq n/m$ this implies that
$$
\left(1-\frac{1}{r_k}\right)n < \frac{n}{m}+2\gamma n+3 \eta n + d_i\left(\frac{n}{m}+3\right)
$$
and therefore, using $3d_i m < 3m^2 \leq n$, $m \geq 1/\gamma$ and $\eta > \gamma$
$$
d_i > \left(1-\frac{1}{r_k}-2\gamma-3\eta\right)m-1-\frac{3d_im}{n} > \left(1-\frac{1}{r_k}-2\gamma-3\eta\right)m-2 \geq
$$
$$
\left(1-\frac{1}{r_k}-4\gamma-3\eta\right)m > \left(1-\frac{1}{r_k}-7\eta\right)m.
$$
This completes the proof of the claim.

We now have that all but at most $\gamma m < 7\eta m$ vertices of $C$ have degree
at least $(1-1/r_k-7\eta)m$ in $C$. Thus, by Corollary \ref{c24}, with $\beta=7\eta$
we have that $C$ contains at least $m/(tr_k)-42\eta r_k m$ vertex disjoint copies of $K(t,r_k)$.
Notice that we can use Corollary \ref{c24} since $m \geq 1/\gamma \geq T(t,r_k,7\eta)$.

Fix a set $L^*$ of at least $m/(tr_k)-42\eta r_k m$ vertex disjoint copies of $K(t,r_k)$ in $C$.
Now, orient the edges of $C$ as follows. The edge $ij \in C$ corresponds
to the fact that $(V_i,V_j)$ is a $\gamma$-regular pair and either
$d(V_i,V_j) \geq \eta$ or $d(V_j,V_i) \geq \eta$. Orient $ij$ from $i$ to $j$ if
$d(V_i,V_j) \geq \eta$, else orient it from $j$ to $i$. In the case that the densities in both directions
are at least $\eta$ we orient the edge arbitrarily. Let ${\vec C}$ denote the resulting orientation of $C$.
${\vec C}$ induces an orientation ${\vec S}$ of each $S \in L^*$.
Since ${\vec S}$ is an orientation of $K(t,r_k)$ we have, by Lemma \ref{l21}, that
${\vec S}$ contains vertex-disjoint copies of $TT_k$ covering all but at most
$f^*(k)-1$ vertices of ${\vec S}$. Over all, we have that ${\vec C}$ contains vertex disjoint
copies of $TT_k$ that cover all vertices of $C$ but at most $42\eta r_k m + (f^*(k)-1)m/(tr_k)$.

Fix a set $L$ of copies of $TT_k$ in ${\vec C}$ that cover all vertices of $C$ but at most $42\eta r_k m + (f^*(k)-1)m/(tr_k)$.
Consider some $S \in L$, and assume, without loss of generality, that the vertices of $S$ are $\{1,\ldots,k\}$
and that each edge of $S$ is directed from a lower vertex to a higher one.
We now show that the subgraph of ${\vec G}$ induced on $V_1 \cup \cdots \cup V_k$ contains
vertex disjoint copies of $TT(h,k)$ that cover all but at most $\epsilon |V_i|/2$ vertices from each
$V_i$, $i=1,\ldots,k$. Proofs of the same nature often appear in applications of
the regularity lemma (see, e.g., \cite{KoSi}). We separate the proof into two lemmas.

\begin{lemma}
\label{l32}
Let $W_1, \ldots ,W_k$ be subsets of vertices having the same size $w$.
Assume that for $i < j$ all edges between $W_i$ and $W_j$ are oriented from
$W_i$ to $W_j$ and that $(W_i,W_j)$ is a $\mu$-regular pair
with $d(W_i,W_j) \geq \eta/2$. If
$$
(k-1) \mu +\frac{h-1}{w}< \left(\frac{ \eta}{4}\right)^{hk}
$$
then there is a $TT(h,k)$ whose color classes $A_1, \ldots, A_k$ satisfy
$A_i \subset W_i$ for $1 \leq i \leq k$.
\end{lemma}
{\bf Proof}\,
We prove that for every $p$, $1 \leq p \leq k$, and for every $q$,
$0 \leq q \leq h$, there are (possibly empty) subsets
$A_i \subset B_i \subset W_i$, $(1 \leq i \leq k)$, 
with the following properties.

\noindent
(i) $|A_i|=h$ for all $i <p$, 
$|A_p|=q$ and $|A_i|=0$ for all $i>p$.

\noindent
(ii) $|B_i| \geq (\frac{\eta}{4})^{(i-1)h}w$ for all $1\leq i \leq p$
and $|B_i| \geq (\frac{\eta}{4})^{(p-1)h+q}w$ for all $p<i \leq k$.

\noindent
(iii)  For all $1 \leq i <j \leq k$, every vertex $u \in A_i$
has an outgoing edge towards every vertex $v \in B_j$.

The assertion of the lemma follows from the above statement
for $p=k$ and $q=h$, since for these values of the parameters the
sets $A_i$ are the color classes of the required $TT(h,k)$.

The subsets $A_i$ and $B_i$ are constructed 
by induction on $(p-1)h+q$.
For $p=1$ and $q=0$ simply
take $A_i=\emptyset$ and $B_i=W_i$ for all $i$.
Given the sets $A_i$, $B_i$ satisfying (i), (ii) and (iii) for 
$p$ and $q$ we show how to modify them for the next value
of $(p-1)h+q$. If $q=h$ and $p<k$ we can replace $p$ by $p+1$
and $q$ by $0$ with no change in the sets $A_i$, $B_i$. Thus we
may assume that $q$ is strictly smaller than $h$.  Consider the set
$D_p=B_p \setminus A_p$. 
Observe that by assumption the size of each 
$B_j$, for $p<j \leq k$ is bigger than $\mu w$.
For each such $j$, let $D_p^j$ denote the set of all vertices
in $D_p$ that have less than $(\eta/2-\mu)|B_j|$ neighbors
in $B_j$. We claim that 
$|D_p^j| < \mu w$ for each $j$. This is because
otherwise the two sets $X=D_p^j$ and $Y=B_j$ would contradict the
$\mu$-regularity of the pair $(W_p,W_j)$, since 
$d(D_p^j,B_j)<\eta/2-\mu$, whereas $d(W_p,W_j) \geq \eta/2$,
by assumption. Therefore,
the size of the set $D_p \setminus (D_p^{p+1} \cup \cdots 
\cup D_p^{k})$
is at least 
$$
|B_p | -|A_p|-(k-p) \mu w \geq 
\left(\frac{\eta}{4}\right)^{(p-1)h}w-q-(k-1)\mu w >0,
$$
where the last inequality follows from the assumption in the lemma.
We can now choose arbitrarily a vertex $v$ in
$D_p \setminus (D_p^{p+1} \cup \cdots \cup D_p^k)$, add it to
$A_p$, and replace each $B_j$ for $p<j \leq k$ by the 
set of neighbors of $v$ in $B_j$. Since 
$\eta/2-\mu >\eta /4$ this will not decrease the
size of each $B_j$ by more than a factor of
$\eta/4$ and it is easily seen that the new sets $A_i$,
$B_i$ defined in this manner satisfy the conditions 
(i), (ii) and (iii) with $p'=p$ and $q'=q+1$. \npf

\begin{lemma}
\label{l33}
Let $V_1, \ldots ,V_k$ be subsets of vertices each of size $b$ or $b+1$.
Assume that for $i < j$ all edges between $V_i$ and $V_j$ are oriented from
$V_i$ to $V_j$ and that $(V_i,V_j)$ is a $\mu^2$-regular pair
with $d(V_i,V_j) \geq \mu+\eta/2$. If
$$
(k-1) \mu +\frac{h-1}{\mu b}< \left(\frac{ \eta}{4}\right)^{hk}
$$
then the graph induced by $V_1 \cup \cdots \cup V_k$
contains at least $(1-2\mu)b/h$ vertex disjoint copies of
$TT(h,k)$, each having $h$ vertices in each $V_i$.
\end{lemma}
{\bf Proof}\,
Let $F$ be a maximal family of vertex disjoint copies of $TT(h,k)$
each having $h$ vertices in each $V_i$.
We prove that the size of $F$ is at least $(1- 2\mu) b/h$. Suppose this  
is false. Let $W_i$ be a subset of $b-h|F|$ vertices of $V_i$ not
appearing in any member of $F$.
Notice that $|W_i|=b-h|F| \geq 2\mu b > \mu(b+1) \geq \mu|V_i|$.
We claim that for all $1 \leq i < j \leq k$, the pair $(W_i,W_j)$ is $\mu$-regular.
First notice that $d(W_j,W_i)=0$ as all edges go from $W_i$ to
$W_j$. Next, notice that every $X \subset W_i$ satisfying
$|X| \geq \mu |W_i|$ also satisfies $|X| \geq \mu^2 |V_i|$, and similarly
every $Y \subset W_j$ satisfying
$|Y| \geq \mu |W_j|$ also satisfies $|Y| \geq \mu^2 |V_j|$.
Therefore
$$
|d(X,Y) - d(W_i,W_j)| \leq |d(X,Y)-d(V_i,V_j)|+|d(V_i,V_j)-d(W_i,W_j)| \leq \mu^2+\mu^2 < \mu.
$$
We have shown the pair $(W_i,W_j)$ is $\mu$-regular.
Notice also that $d(W_i,W_j) \geq \eta/2$.
By Lemma \ref{l32}, with $w=b-h|F|$,
there must be another copy of $TT(h,k)$, disjoint from $F$,
contradicting its maximality. \npf

We apply Lemma \ref{l33} to the subgraph of ${\vec G}$ induced on $V_1 \cup \cdots \cup V_k$
(ignoring directed edges going in the ``wrong'' direction, from a larger indexed class to a smaller one).
We are allowed to do this since $\mu^2 \geq \gamma$, $d(V_i,V_j) \geq \eta > \eta/2+\mu$,
and
$$
(k-1) \mu +\frac{h-1}{\mu b}= (k-1) \mu +\frac{h-1}{\mu \lfloor n/m \rfloor}< (k-1) \mu + (h-1) \mu <
 \left(\frac{ \eta}{4}\right)^{hk}.
$$
By Lemma \ref{l33} we obtain vertex disjoint copies of $TT(h,k)$ that cover all vertices of $V_i$ but at most
$2\mu b + 1 < \epsilon b/2$. Repeating this process for every $S \in L$ we obtain vertex
disjoint copies of $TT(h,k)$ that cover all the $n$ vertices of $G$ but at most
$\epsilon b/2 \cdot m + (b+1)(42\eta r_k m + (f^*(k)-1)m/(tr_k))$.
Since
$$
\epsilon \frac{b}{2} \cdot m + (b+1)\left(42\eta r_k m + (f^*(k)-1)\frac{m}{tr_k}\right)
\leq \frac{\epsilon}{2}n + 2n\left(42\eta r_k + \frac{2^{k-1}}{tr_k}\right) < \epsilon n
$$
the theorem follows. \npf

\section{Lower bounds}
{\bf Proof of Proposition \ref{p2}:}\,
Let $1/30 > \gamma > 0$. We show there exists a graph $G$ with minimum degree at least $(5/6-\gamma)n$ and
an orientation of $G$ having at most $n/3-\gamma n$ vertex-disjoint copies of $TT_3$. Clearly we may assume
that $\gamma$ is rational. Let $\alpha=\gamma+1/6$ and let $n$ be chosen such that
$\alpha n$ is an integer. Let $G$ be the complete $6$-partite graph with $n$ vertices and with vertex partition
$V_1,\ldots,V_6$ where $|V_i|=\alpha n$ for $i=2,\ldots,6$ and $|V_1|=(1-5\alpha)n$. Notice that since
$\alpha < 1/5$ we have that $V_1 \neq \emptyset$. Also notice that $\delta(G)=(1-\alpha)n=(5/6-\gamma)n$.
Let $T$ be a tournament on 6 vertices, where $(5,6)$ is an edge of $T$, and for each $i=1,2,3,4$, both $(6,i)$ and $(i,5)$
are edges of $T$. The orientation of the other 6 edges of $T$ can be chosen arbitrarily.
Now consider the orientation ${\vec G}$ of $G$ where all edges between $V_i$ and $V_j$ are directed from $V_i$
to $V_j$ if and only if $(i,j)$ is an edge of $T$. By construction, any copy of $TT_3$ in ${\vec G}$ has
at most one vertex in $V_5 \cup V_6$. Thus, there are at most
$(|V_1|+|V_2|+|V_3|+|V_4|)/2=(1/2-\alpha)n=n/3-\gamma n$
vertex-disjoint copies of $TT_3$ in ${\vec G}$. \npf

We now show that for $k \geq 4$, we cannot replace the constant $c_k=1-2^{-k-\log k}$ appearing in Theorem
\ref{t1} with a constant less than $1-2^{-0.5k(1+o(1))}$. Recall the result of Erd\H{o}s and Moser, mentioned
in the introduction, stating that $f(k) \geq 2^{0.5k(1+o(1))}$. Let $T$ be a tournament with $f(k)-1$
vertices, not containing $TT_k$ as a subgraph. Consider the complete $(f(k)-1)$-partite graph $G$ with
$n/(f(k)-1)$ vertices in each part. The minimum degree of $G$ is $n(1-1/(f(k)-1)) \geq n(1-2^{-0.5k(1+o(1))})$.
Consider the orientation ${\vec G}$ of $G$ formed by replacing each vertex of $T$ with
an independent set of size $n/(f(k)-1)$. Clearly, ${\vec G}$ does not have even a single copy of
$TT_k$ as a subgraph. In fact, Erd\H{o}s and Moser conjectured that $f(k) \geq 2^{k(1-o(1))}$. If this conjecture
is true then the constant $c_k$ in Theorem \ref{t1} is rather tight.

\section{Exact factors}
{\bf Proof of Theorem \ref{t2}:}\,
Let $G$ have $n$ vertices, $n=ht$, and $\delta(G) \geq n(1-1/4^h)+4^h$.
Let ${\vec G}$ be an orientation of $G$.
Recalling the definition of $g(k)$ from the introduction, and the fact that $g(k) < 4^k$, we have,
By Tur\'an's Theorem (cf. \cite{Bo}), that $G$ contains a complete graph on $g(h)$ vertices.
Thus, ${\vec G}$ contains a tournament on $g(h)$ vertices. By definition, this tournament
has a $TT_h$-factor, and, in particular, an $H$-factor with $g(h)/h$ copies of $H$.
Let $m \equiv n \bmod g(h)$. Since $h | n$ and $h | g(h)$ we have $h |m$. Pick $m/h$
vertex disjoint copies of $H$ from the obtained $H$-factor and delete their vertices from $G$.
We remain with a graph $G'$ (and its corresponding orientation ${\vec G'}$) on $n'=n-m$ vertices,
where $n' \equiv 0 \bmod g(h)$. Furthermore,
$$
\delta(G') \geq \delta(G) - m \geq n\left(1-\frac{1}{4^h}\right) \geq n'\left(1-\frac{1}{4^h}\right) \geq n'\left(1-\frac{1}{g(h)}\right).
$$
By the theorem of Hajnal and Szemer\'edi, $G'$ has a $K_{g(h)}$-factor.
In particular ${\vec G'}$ contains $n'/g(h)$ tournaments, each having $g(h)$ vertices.
By definition of $g(h)$, each of these tournaments has a $TT_h$-factor and hence also an $H$-factor. \npf

For $H=TT_2$ and $H=TT_3$ it is very easy to determine a sharp analog of Theorem \ref{t2}.
Every graph with $n$ vertices, $n$ even, and minimum degree $n/2$ has a perfect matching.
Thus, in every orientation there is a $TT_2$-factor. If $n \equiv 0 \bmod 3$ and the minimum degree
of an $n$-vertex graph $G$ is $5n/6$ then there are two cases. If $n \equiv 0 \bmod 6$ then by the
Hajnal and Szemer\'edi Theorem, $G$ has a $K_6$-factor. Since $g(3)=6$ we have that
every orientation of $G$ has a $TT_3$-factor. If $n \equiv 3 \bmod 6$ then we proceed as follows.
If ${\vec G}$ is an orientation of $G$ then we pick one (of course there is one) $TT_3$ and delete its vertices.
We now have $n-3$ vertices and minimum degree at least $\lceil 5n/6 \rceil-3 \geq 5(n-3)/6$.
As in the previous case, there is a $K_6$-factor in the undirected remaining graph
and a $TT_3$-factor in the directed remaining graph. The proof of Proposition \ref{p2}
shows that the constant $5/6$ is optimal.

\section{Concluding remarks and open problems}
\begin{itemize}
\item
Theorem \ref{t1} is asymptotically optimal for acyclic digraphs $H$ with ${\vec \chi(H)} \leq 3$. It would be interesting
to determine a sharp minimum degree requirement also for acyclic digraphs with ${\vec \chi(H)} \geq 4$.
Notice that the proof of Theorem \ref{t1} shows that it suffices to prove a sharp minimum degree requirement
for $TT_k$, and the same bound would hold for all fixed graph $H$ with ${\vec \chi(H)} = k$.
Even for $TT_4$ this is still open. The comment after lemma \ref{l21} and the proof of Theorem \ref{t1} yield an
upper bound of $19/20$ for  $TT_4$. Namely, for every $\epsilon > 0$, if $n$ is sufficiently large and
$G$ has $n$ vertices and minimum degree at least $19n/20$ then every orientation of $G$ has vertex disjoint copies
of $TT_4$ covering all but at most $\epsilon n$ vertices. The best lower bound that we currently have is $25/28$.
This is obtained as follows. Let $T$ be a tournament on 7 vertices without a $TT_4$ (in fact $T$ is
unique \cite{RePa}). Let $G$ be the 10 partite graph $G$ with $n$ vertices and with vertex classes $V_1, \ldots, V_{10}$.
where $|V_i|=(1+\gamma)3n/28$ for $i=1,\ldots,7$, $|V_i|=(1-3\gamma)n/12$ for $i=8,9,10$.
Notice that $\delta(G) =(25-3\gamma)n/28$.
Consider the orientation ${\vec G}$ of $G$ where for $1 \leq i < j \leq 7$, the edges between $V_i$ and $V_j$
all go in the same direction, corresponding to the direction of the edge $ij$ of $T$. All other edges are oriented arbitrarily.
Each $TT_4$ of $G$ must contain at least one vertex from $V_8 \cup V_9 \cup V_{10}$.
Thus, the maximum number of vertex-disjoint $TT_4$ in ${\vec G}$ is at most $(1-3\gamma)n/4$. Hence,
at least $3\gamma n$ vertices are uncovered.
\item
In the undirected case, the minimum degree guaranteeing an almost $K_k$-factor is
the same as the minimum degree guaranteeing an exact $K_k$-factor. In fact, the
$K_k$-factor problem is completely settled in the Hajnal and Szemer\'edi theorem,
and this theorem does not require the use of Szemer\'edi's regularity lemma, and works for all $n$.
This is also true in the directed case for $TT_2$ and $TT_3$, as shown in Section 5. However, we currently have no
proof yielding the bound $c_k=1-2^{-k -\log k}$ of Theorem \ref{t1} that avoids the use of the
(directed) regularity lemma even for $TT_4$. In fact, we conjecture that for $k \geq 4$ there is a gap
between the threshold guaranteeing a $TT_k$-factor and the threshold guaranteeing an {\em almost} $TT_k$-factor.
\item
It would be interesting to determine an asymptotically tight minimum degree requirement for an almost $H$-factor
for every acyclic digraph $H$, in terms of its chromatic number and the sizes of its vertex classes.
Currently, theorem \ref{t1} is only applicable to acyclic digraphs with balanced $k$-colorings and is
tight for $k=2,3$. In fact, even providing an analog of theorem \ref{t1} applicable to all bipartite acyclic digraphs
would be interesting.
\end{itemize}

\end{document}